\documentclass{amsart}
\usepackage{amssymb}
\usepackage[colorlinks=false,breaklinks=true,pdfpagemode=UseThumbs,
bookmarksnumbered=true]{hyperref}

\copyrightinfo{}{}

\newtheorem{theorem}{Theorem}[section]
\newtheorem{lemma}[theorem]{Lemma}

\theoremstyle{definition}
\newtheorem{definition}[theorem]{Definition}
\newtheorem{example}[theorem]{Example}

\theoremstyle{remark}
\newtheorem{remark}[theorem]{Remark}

\numberwithin{equation}{section}

\begin{document}

\title[Equivalence Theorems in Numerical Analysis]
{Equivalence Theorems in Numerical Analysis :\\ 
Integration, Differentiation and Interpolation}
\author{John Jossey}
\address{
  Department of Computer Science\\
  University of Illinois at Urbana-Champaign\\
  201 N. Goodwin Ave., Urbana, IL 61801}
\email{johnjossey@gmail.com}

\author{Anil N. Hirani}
\address{
  Department of Computer Science\\
  University of Illinois at Urbana-Champaign\\
  201 N. Goodwin Ave., Urbana, IL 61801}
\email{hirani@cs.uiuc.edu}
\urladdr{http://www.cs.uiuc.edu/hirani}
\thanks{Author for correspondence : Anil N. Hirani (hirani@cs.uiuc.edu).}

\subjclass[2000]{Primary 65J05, 65J10, 65J15}
\keywords{Consistency, stability, convergence, equivalence theorems.}
\date{}

\begin{abstract} 
  We show that if a numerical method is posed as a sequence of
  operators acting on data and depending on a parameter, typically a
  measure of the size of discretization, then consistency, convergence
  and stability can be related by a Lax-Richtmyer type equivalence
  theorem -- a consistent method is convergent if and only if it is
  stable. We define consistency as convergence on a dense subspace and
  stability as discrete well-posedness. In some applications
  convergence is harder to prove than consistency or stability since
  convergence requires knowledge of the solution. An equivalence theorem
  can be useful in such settings. We give concrete instances of
  equivalence theorems for polynomial interpolation, numerical
  differentiation, numerical integration using quadrature rules and
  Monte Carlo integration.
\end{abstract}

\maketitle

\section{Introduction}
For a numerical method the three most important aspects are its
consistency, convergence and stability. These three were related in
the well known equivalence theorem of Lax and Richtmyer for finite
difference methods for certain partial differential equations
\cite{LaRi1956}. We show that in a very general setting of numerical
methods, in which a numerical method is posed as a family of operators
acting on data, there is a Lax-Richtmyer type equivalence theorem : a
consistent method is convergent if and only if it is stable.  After
proving the theorem in two general settings
(Theorem~\ref{thm:equiv_linear} and
Theorem~\ref{thm:equiv_nonlinear}), we prove it in specific instances
of numerical integration using quadrature, numerical integration using
Monte Carlo methods, numerical differentiation and polynomial
interpolation.

Consistency is a measure of how good the discretization is.  Roughly,
it says that the discretization is close to the smooth operator in
some sense. If the discrete solution converges to the smooth solution
then the numerical method is said to be convergent. Note that for
discussing consistency and convergence one needs some information
about the smooth problem and the smooth solution. However, numerical
stability is purely a property of the discrete scheme. Roughly,
stability means that the propagated error is controlled by the error
in the data. Hence there is a similarity between numerical stability
and well-posedness.

In practice, convergence can be the hardest to prove among
consistency, convergence and stability, since the actual solution is
usually not known. Hence equivalence theorems can be useful in such
situations.  Equivalence theorems essentially say that we need not
worry about the convergence while solving a problem numerically as
long as its discretization is consistent with the smooth problem and
the discrete scheme is stable. In addition, such theorems also show
that unstable schemes will not converge for some data. These are the
same advantages that are often appreciated in the setting of the
classical Lax-Richtmyer equivalence theorem for finite difference
schemes (see \cite{Strikwerda2004}, page 32).

After the preliminaries in the next section, in Section~\ref{sec:ccs}
we define consistency, convergence and stability in a general context
of operators acting on data and we then prove equivalence theorems in
this setting. The three sections that follow specialize these notions
to specific classes of basic numerical methods. The convergence and
stability theory for these example areas are well understood, but
equivalence theorems have not been discussed in these areas in the
literature.

As given here, the equivalence theorems in these example areas serve
only as concrete instances for illustration of the main ideas. We make
no claims that these examples have direct practical importance in
numerical analysis. However, with proper generalizations, the ideas
might be of use in practical situations. For example, convergence of
multidimensional interpolation can be related to its stability and
consistency, as we sketch in Section~\ref{sec:interp}.  Until now the
advantages of equivalence theorem have been limited to finite
difference methods. We suggest that similar benefits may be possible
in many areas of numerical analysis.

\section{ Preliminaries}
\label{sec:prelim}
For the convenience of the reader, we state some definitions and
theorems used later on in the paper.  Uniform boundedness principle is
one of the fundamental building blocks of functional analysis and it
is useful for proving equivalence theorems in the linear operator
setting. It says that a sequence of pointwise bounded continuous
linear operators defined on a complete normed linear space are
uniformly bounded.

\begin{theorem}[\bf Uniform Boundedness Principle]
  Let $\{F_{i \in I}\}$ be a set of bounded linear operators from a
  Banach space $V$ to a normed linear space $W$, where $I $ is an
  arbitrary set.  Assume for every $v \in V$, the set $\{F_{i }(v)\}$
  is bounded. Then $\sup_{i \in I} \left \Vert F_i \right \Vert <
  \infty$.
\end{theorem}

\begin{proof} 
  See \cite{AbMaRa1988} or \cite{Conway1990}.  The main ingredient of
  the proof is Baire category theorem.
\end{proof}

The next lemma will be used in Section~\ref{sec:interp} for proving that 
polynomial interpolation operators are bounded.

\begin{lemma}\label{lem:closed_kernel}
  Let $V$ and $W$ be normed linear spaces over $\mathbb{R}$, where $W$
  is finite dimensional.  Let $T : V \rightarrow W$ be a surjective
  linear operator with a closed kernel.  Then $T$ is continuous.
\end{lemma}
\begin{proof}
  Since $K = T^{-1}(0)$ is closed, $V/K$ is a normed linear space (see
  for instance Theorem 4.2 on page 70 in \cite{Conway1990}).  For $v
  \in V$ the norm in $V/K$ is defined as usual to be $\Vert v + K
  \Vert_{V/K} := \inf\{\Vert v + k \Vert\,\text{s.t.}\, k \in K\}$,
  which is the distance of $v$ from $K$.  Let $\widehat{T}: V/K
  \rightarrow W$ be the unique linear map such that $T = \widehat{T}
  \circ \pi$, where $\pi: V \rightarrow V/K$ is the quotient
  map. Recall that the quotient map is continuous. Note also that
  $\widehat{T}$ is a bijection. Then since $\widehat{T}$ is a linear
  bijection between $V/K$ and $W$, and $W$ is finite dimensional, we
  have that $V/K$ is finite dimensional. Thus $\widehat{T}$ is
  continuous (see page 56 of \cite{AbMaRa1988}). Hence $T$ is
  continuous since it is the composition of two continuous functions.
\end{proof}

\begin{remark}
  Note that all that one needs above is that the image of $T$ in $W$
  be finite dimensional and the kernel be closed. The surjectivity of
  $T$ is not required in that case.
\end{remark}

The remaining part of this section deals with some basic facts about
probability theory which are needed when we discuss Monte Carlo
integration. These are not needed for the general equivalence theorems
of Section~\ref{sec:ccs} or other sections with the exception of
Section~\ref{subsec:mc}. Let $(\Omega, \Sigma, P)$ denote a
probability space \cite{Varadhan2001} where $\Omega$ is the space of
outcomes, $\Sigma$ is a $\sigma$-algebra of subsets of $\Omega$, and
$P$ is a countably additive probability measure on $\Sigma$.

\begin{definition}
  A random variable $X: \Omega \rightarrow \mathbb{R}$ is a measurable
  function, i.e., for every Borel set $B \subset \mathbb{R}$,
  $X^{-1}(B) \in \Sigma$.
\end{definition}

The probability measure $\alpha = PX^{-1}$ defined on $\mathbb{R}$ is
called the distribution of $X$.  We will assume that the random
variable $X$ is continuous, i.e., there exists a nonnegative function
$f(x)$, called the probability density function (pdf) of the random
variable $X$, defined on $\mathbb{R}$ such that for all Borel sets $A
\subset \mathbb{R}$
\[ 
P\left[\omega : X(\omega) \in A \right]=\alpha(A) = \int_{A}f\, .
\]
The mean or the expectation of the random variable $X$ if it exists is
$E(X) = \int_{\mathbb{R}} x f(x) dx$.

\begin{definition}
  Two random variables $X, Y$ are said to be independent if
  \[
  P[\omega : X(\omega) \in A, Y(\omega) \in B] =P[\omega : X(\omega)
  \in A]P[\omega : Y(\omega) \in B]\, ,
  \]
  for all Borel sets $A,B$.
\end{definition}

\begin{theorem}[\bf Strong Law of Large Numbers]\label{thm:lawof}
  Let $ X_1, X_2, \ldots, X_n, \ldots$ be a sequence of independent
  and identically distributed random variables with finite mean
  $\mu$. Then
  \[
  P\left[\omega : \lim_{n \rightarrow \infty}\dfrac{1}{n}\sum
    _{i=1}^{n}X_i(\omega) =\mu \right]=1\, .
  \]
\end{theorem}

\begin{proof}
  See \cite{Chung2001}.
\end{proof}

\begin{theorem} \label{thm:independent} Let $X$,$Y$ be two
  independently and identically distributed random variables. If $f$
  is a Borel measurable function \cite{Rudin1987} on $\mathbb{R}$,
  then the random variables $f(X)$,$f(Y)$ are independent and
  identically distributed.
\end{theorem}

\begin{proof}
  Let $A, B$ be any Borel sets in $\mathbb{R}$.
  \begin{align*}
    P\left[\omega : f(X(\omega)) \in A, f(Y(\omega)) \in B \right]
    &=P\left[\omega : X(\omega) \in f^{-1}(A), Y(\omega)
      \in f^{-1}(B) \right]\\
    &=P\left[\omega : X(\omega) \in f^{-1}(A) \right]
    P\left[\omega : Y(\omega) \in f^{-1}(B) \right]\\
    &=P\left[\omega : f(X(\omega)) \in A \right] P\left[\omega :
      f(Y(\omega)) \in B \right]\, .
  \end{align*}
  Hence $f(X), f(Y)$ are independent. 

  Let $\alpha $ be the distribution measure of both $X$ and $Y$, i.e.,
  \[\alpha(A) = P \left[\omega: X(\omega) \in A \right]
  = P \left[\omega: Y(\omega) \in A \right]\, .\]Let $\beta_{1}$ $(
  \beta_{2})$ be the distribution of $f(X)$ $(f(Y)$ respectively).
  Therefore,
  \[\beta_{1}(A)=P\left[\omega : f(X(\omega)) \in A \right]=
  P\left[\omega : X(\omega) \in f^{-1}(A) \right]=\alpha(f^{-1}(A))\,
  ,\]
  \[\beta_{2}(A)=P\left[\omega : f(Y(\omega)) \in A \right]=
  P\left[\omega : Y(\omega) \in f^{-1}(A) \right]=\alpha(f^{-1}(A))\,
  .\] Hence $\beta_1 = \beta_2 (= \alpha f^{-1})$.
\end{proof}

\begin{remark}
  We need the condition of Borel measurability on $f$ because for
  every Borel set $A$, we want $f^{-1}A$ to be a Borel set so that
  $X^{-1}(f^{-1}A) \subseteq \Sigma$. Otherwise, if $f$ is not a Borel
  measurable function, then $f^{-1}A$ need not be a Borel measurable
  set and then $X^{-1}(f^{-1}A)$ would not be in $\Sigma$.
\end{remark}

\section{Consistency, Convergence and Stability}
\label{sec:ccs}

If a smooth problem can be formulated as an operator applied to data,
then the discretization can usually be formulated as a family of
operators depending on some discretization parameter, typically a
measure of the mesh size. For example, the parameter might be a
measure of the distance between nodes in quadrature, or between
interpolation points in polynomial interpolation, etc. 

The smooth and discrete operators can be made to act on the same
space. In some cases this may be done for example, by considering
continuous functions instead of discrete data. We then define the
discrete scheme to be convergent if the discrete operators converge to
the smooth operator pointwise on the entire space. We define
consistency to be convergence on a dense subspace. If the discrete
operators are bounded linear then the definition of stability is
uniform boundedness of the family of discrete operators. However, when
the discrete operators are general nonlinear operators, we define
stability as asymptotic pointwise boundedness of the family of
operators. The precise definitions are given below in
Definition~\ref{defn:ccs_linear} and~\ref{defn:ccs_nonlinear}.

The two definitions of stability lead to two different proofs for the
equivalence theorem, both of which appear in this section.
Theorem~\ref{thm:equiv_linear} is a Lax-Richtmyer type equivalence
theorem applicable to general numerical analysis problems when the
discrete operators involved are bounded linear.  When this condition
is dropped, we get Theorem~\ref{thm:equiv_nonlinear}. Specific
incarnations of the linear case theorem are proved later in the
context of numerical integration (Theorems~\ref{thm:quadrature}),
numerical differentiation (Theorem~\ref{thm:differentiation}), and
polynomial interpolation (Theorem~\ref{thm:interp}). We treat the
general problem of Monte Carlo integration without the assumption of
linearity. This leads to a proof of an equivalence theorem
(Theorem~\ref{thm:montecarlo}) analogous to the nonlinear case but
with a probabilistic flavor.

In the classical Lax-Richtmyer equivalence theorem for partial
differential equations, the main ingredient in the proof of stability
implying convergence is essentially triangle inequality and a density
argument \cite{LaRi1956}.  Similarly, to prove the equivalence
theorems in this paper, the main tools we use are uniform boundedness
principle, triangle inequality and some density arguments. 

Let $V$ be a Banach space and $W$ a normed linear space.  Let $h \in
(0,1)$. The upper limit of $h$ is not relevant because we will be
considering limits as $h \rightarrow 0$. Let $T,T_{h} : V \rightarrow
W$ be set of bounded linear operators.

\begin{definition}[Linear Discrete Operator Case]\label{defn:ccs_linear}
  If $\lim_{h \rightarrow 0} \Vert(T_h-T)v\Vert =0$ for every $v \in
  V$, then $T_h$ is said to {\bf converge} to $T$, and if $\lim_{h
    \rightarrow 0} \Vert(T_h-T)v\Vert =0$ for every $v \in V_0$, where
  $V_0$ is a dense subspace of $V$, then $T_h$ is said to be {\bf
    consistent} with $T$. If $\sup_{h} \Vert T_h \Vert < \infty$ then
  $T_h$ is called {\bf stable}.
\end{definition}

Our definition of consistency is motivated by the definition of
consistency for finite difference schemes for certain PDEs. For a
partial differential equation $P u = f$, where $u$ is the unknown, and
a finite difference scheme $P_{k,h} v = f$, the finite difference
scheme is called consistent if for any smooth function $\phi(t,x)$,
$P\phi - P_{k,h}\phi \rightarrow 0$ as $k,h \rightarrow 0$. Here $k$
and $h$ are a measure of the space and time meshes respectively. The
convergence is pointwise convergence at every $(t,x)$. This definition
is from \cite{Strikwerda2004}. Although density is not explicitly
mentioned in this definition, note that smooth functions are dense in
the typical function spaces in which the solutions live.

\begin{remark}\label{rem:consistency}
  Observe that by our definition of consistency, convergence implies
  consistency. However, there are other definitions of consistency
  under which there exist inconsistent schemes which converge. See for
  instance, \cite{Yamamoto2002} and Example 1.4.3 in
  \cite{Strikwerda2004} in the context of finite difference schemes
  for PDEs.
\end{remark}

The stability definition above is equivalent to discrete
well-posedness, i.e., well-posedness of the discrete problems which is
that for any $h$ and for any $v_1,v_2$ in $V$, $\Vert T_h(v_1 - v_2)
\Vert \le K \Vert v_1 - v_2 \Vert$ where $K \ge 0$ is some constant.
 
\begin{theorem}[\bf Equivalence Theorem for Linear Discrete Operators]
  \label{thm:equiv_linear}
  A consistent family of operators $T_h$ is convergent if and only if
  it is stable.
\end{theorem}
\begin{proof}
  Suppose the family is convergent, i.e., $\lim_{h \rightarrow 0} T_h
  v =Tv$ for every $v \in V$. Hence $\Vert T_h v\Vert \leq K( v)$
  where $K(v)$ is a constant possibly depending on $v$.  Since $V$ is
  complete, by uniform boundedness principle we have uniform
  boundedness of $T_h$, i.e., $\sup_{h } \Vert T_h \Vert\leq K $, for
  some $K \geq 0$, and hence stability.

  Conversely, suppose that $T_h$ is consistent and stable. Since $V_0$
  is a dense subspace of $V$, for a given $v \in V$ choose $v_0 \in
  V_0$ such that
  \[
  \Vert v- v_0 \Vert \leq \dfrac{\epsilon}{3 \max \{\left\Vert T
    \right \Vert, \sup_{h} \left\Vert T_h \right\Vert\}}\, .
  \]
  Because of consistency, there exists $ h_0 \in (0,1)$ such that, for
  all $h \leq h_0$ we have that $\Vert T_h v_0 - T v_0 \Vert \leq
  \epsilon /3$. Hence for all $h \leq h_0$
  \begin{align*}
    \Vert T_h v - T v \Vert &\leq \Vert T_h v - T_h v_0\Vert +
    \Vert T_h v_0 - T v_0 \Vert+\Vert T v_0 -  T v \Vert \\
    & \leq \Vert T_h \Vert \Vert v-v_0\Vert + \Vert T_h v_0 - T v_0
    \Vert+
    \Vert T \Vert \Vert v_0 - v\Vert \\
    & \leq \dfrac{\epsilon}{3}+\dfrac{\epsilon}{3}+\dfrac{\epsilon}{3}
    = \epsilon \, .
  \end{align*}
  Therefore we have convergence.
\end{proof}

For the case when the discrete operators are nonlinear we give a
different definition of stability since uniform boundedness is
unlikely. This allows us to prove an equivalence theorem in this case
as well. Most of the examples we give in later sections are linear
ones.  The proof of the equivalence theorem in the case of nonlinear
discrete operators (Theorem~\ref{thm:equiv_nonlinear}) is different
from the linear case.

Let $T_n$ be a sequence of operators (not necessarily linear or
continuous) from $V$ to $W$, where $V$ and $W$ are normed linear. Let
$T: V\rightarrow W$ be a bounded linear operator.

\begin{definition}[Nonlinear Discrete Operator Case]\label{defn:ccs_nonlinear}
  The discrete operator $T_n$ is said to {\bf converge} to $T$ if, for
  any given $v \in V$, the sequence $T_n v$ converges to $Tv$ in $W$,
  i.e., for each $v \in V$, given $\epsilon > 0$ there exists $n_0 \in
  \mathbb{N}$, such that for all $n \geq n_0$, ${\left \Vert T_n v - T
      v \right \Vert } \leq \epsilon$.  If there is a dense subspace
  $V_0$ of $V$ such that for any $v$ in $V_0$, the sequence $T_n v$
  converges to $Tv$ in $W$ then $T_n$ is said to be {\bf consistent}
  with $T$. The operator $T_n$ is {\bf stable at} $v_0 \in V$ if for
  any $\epsilon > 0$ there exists a $\delta > 0$, such that for each
  $v \in V$ with ${\left \Vert v - v_0 \right \Vert } \leq \delta$,
  there exists $ n_0 \in \mathbb{N} $ such that for all $n \geq n_0$,
  ${\left \Vert T_n v - T_n v_0 \right \Vert } \leq \epsilon$.  It is
  {\bf stable} if it is stable at every $v_0 \in V$.
\end{definition}

Remark~\ref{rem:consistency} about convergence implying consistency
that was made earlier about the linear case is also valid for the
nonlinear case covered in Definition~\ref{defn:ccs_nonlinear} since the
consistency and convergence definitions are the same in both cases.

\begin{theorem}[\bf Equivalence Theorem for Nonlinear Discrete Operators]
\label{thm:equiv_nonlinear}
A consistent family of operators $T_n$ is convergent if and only if it
is stable.
\end{theorem}

\begin{proof}
  Suppose $T_n$ is convergent. This implies that given $v, v_0$ in $V$
  and $\epsilon > 0$ there exists $n_v, n_{v_{0}} \in \mathbb{N}$ such
  that ${\left\Vert T_n v -T v \right \Vert} \leq \epsilon/3 $ for all
  $n \geq n_v$ and $\left\Vert T_n v_0 -T v_0 \right \Vert \leq
  \epsilon/3$ for all $n \geq n_{v_0}$.  We need to find an $n_0 \in
  \mathbb{N}$ and a $\delta > 0$ such that for a given $v \in V$ with
  $\left\Vert v -v_0\right\Vert < \delta$, $\left\Vert T_n v -T_n v_0
  \right\Vert \leq \epsilon$ for all $ n \geq n_0$.  But, $\left\Vert T_n
    v -T_n v_0 \right\Vert = \left\Vert T_n v - T v + T v - T v_0 + T
    v_0- T_n v_0 \right\Vert \leq \left\Vert T_n v - T v \right\Vert +
  \left\Vert T v-T v_0 \right\Vert + \left\Vert T v_0- T_n v_0 \right
  \Vert$. Since $T$ is bounded linear operator, if we choose $\delta$
  appropriately, then $\Vert T v -T v_0\Vert \leq \Vert T\Vert \Vert v
  - v_0 \Vert \leq \epsilon/3$. Letting $n_0 = \max (n_v, n_{v_0})$ we
  get $\left\Vert T_n v -T_n v_0 \right\Vert \leq \epsilon$ for all $
  n \geq n_0$. Since $ v_0$ was arbitrary, $T_n$ is stable if it is
  convergent.

  Conversely, suppose we assume stability and consistency.  We'll show
  convergence at $v_0$.  Thus we need to find an $n_0 \in \mathbb{N}$
  such that ${\left\Vert T_n v_0 - Tv_0 \right\Vert} \leq \epsilon$
  for all $n \geq n_0$.

  Stability at $v_0$ means that there exists $ \delta >0$ such that
  for each $v_0' \in V_0$ with ${\left\Vert v_0'- v_0 \right\Vert}
  \leq \delta$ there exists $n_1 \in \mathbb{N}$ with ${\left\Vert T_n
      v_0 - T_n v_0' \right\Vert} \leq \epsilon /3$ for all $n \geq
  n_1$.  Since $T$ is bounded linear operator, if we choose $\delta$
  appropriately we can also make $\left\Vert T v_0' - T v_0
  \right\Vert \leq \epsilon/3$. Choose such a $\delta$ and $v_0' \in
  V_0$. Note that that $\left\Vert T_n v_0 - Tv_0 \right\Vert \leq
  \left\Vert T_n v_0 - T_nv_0' \right\Vert + \left\Vert T_n v_0' - T
    v_0' \right\Vert + \left\Vert T v_0'-T v_0 \right\Vert$.  By the
  choice of $v_0'$ the first and last terms on the right hand side of
  the above inequality are already at most $\epsilon/3$.  By
  consistency, there exists $n_2 \in \mathbb{N}$ such that
  ${\left\Vert T_n v_0'- T v_0' \right\Vert}\leq \epsilon/3$ whenever
  $n \geq n_2$. Choose $n_0 = \max (n_1, n_2)$. Then for all $n \ge
  n_0$ we have $\left\Vert T_n v_0 - T v_0 \right\Vert \leq
  \epsilon$. Hence we have convergence at $v_0$.  Since $v_0 $ was
  arbitrary, we have that stability implies convergence.
\end{proof}

In the linear case above (Definition~\ref{defn:ccs_linear} and
Theorem~\ref{thm:equiv_linear}), we used a real parameter
$h$. Typically this will be a measure of size of discretization, such
as the maximum distance between adjacent nodes in the partition of an
interval. In the nonlinear case (Definition~\ref{defn:ccs_nonlinear}
and Theorem~\ref{thm:equiv_nonlinear}) we chose to use a natural
number $n$ as the parameter. This might stand, for example, for the
number of times sampling is done in Monte Carlo integration. This
change from real $h$ to natural number $n$ was done to give both
flavors of the definitions and proofs. Each can be written using
either $h$ or $n$.

\begin{remark}\label{rem:conv_stable}
  Note that in the proofs of both the equivalence theorems above,
  we did not assume consistency to show that convergence implies
  stability. It would have been redundant anyway, to assume consistency
  when we already have convergence (see Remark~\ref{rem:consistency}).
\end{remark}

\begin{remark}
  In \cite{AtHa2005} (page 67) consistency, convergence and stability
  are defined in a general setting of linear operators. The problem
  setting is the solution of equation $L v = w$, where $L$ is a
  bounded linear operator.  This is discretized as $L_n v_n = w$. Here
  the unknown is $v$ and $v_n$ and $w$ is known. Under their
  definitions they show one side of the equivalence theorem, that a
  consistent method is convergent if it is stable.  Our setting
  however is that of ``direct'' problems. In our case the object being
  approximated discretely is $T f$ which is approximated by $T_h f$
  where $f$ is known data. In \cite{AtHa2005} an inverse of the
  operator $L$ is required. In our case $T$ may go from the space of
  continuous functions to reals (as in Section~\ref{sec:integration})
  so that an inverse may not exist.
\end{remark}

\section{Numerical Integration}
\label{sec:integration}
As the first application of the ideas of Section~\ref{sec:ccs} we now
discuss the notions of consistency, convergence and stability of
numerical integration. We only address definite integrals of
continuous functions on the real line. The two most successful methods
for numerical integration are quadrature rules and Monte Carlo
integration and these are covered below in
Sections~\ref{subsec:quadrature} and ~\ref{subsec:mc}.  The
definitions and proof of theorem for quadrature are identical to the
linear case in Section~\ref{sec:ccs}. For Monte Carlo integration
however, the notions of consistency, convergence and stability need to
be put into a probabilistic setting and the equivalence theorem proof
uses some probabilistic reasoning. Otherwise, the pattern of the proof
follows that of Theorem~\ref{thm:equiv_nonlinear}. In the quadrature
case we apply the theorem to infer convergence of Gaussian quadrature
from a simple proof of its stability and we also discuss composite
trapezoidal rule and the instability of Newton-Cotes quadrature. In
the case of Monte Carlo integration we discuss the Sample Mean method
as an example.

\subsection{Quadrature}\label{subsec:quadrature}
The numerical approximation of definite integrals is often done using
quadrature rules \cite{Heath2002}. Let $V = (C[a,b], {\Vert \cdot
  \Vert}_{\infty})$. For $f \in V$ define $I(f) = \int_{a}^{b}f(x)
dx$, which can be approximated by a sequence of quadratures $I_{n}(f)
= \sum_{i=0}^{n} {w_{i}^{(n)}f(x_{i}^{(n)})}$, where $a \leqslant
x_{0}^{(n)} < x_{1}^{(n)} < \cdots < x_{n}^{(n)} \leqslant b$ is a
partition of $[a,b]$. The points $x_i$ are called nodes. These nodes
are not necessarily equally spaced or progressive. (In progressive
quadrature if the number of nodes is increased from $n_1$ to $n_2$
then only $n_2-n_1$ nodes are new.) The real linear functional $I: V
\rightarrow \mathbb{R}$ is a bounded linear operator and $\left\Vert I
\right\Vert = b-a$.

Each $I_n: V \rightarrow \mathbb{R}$ is also a linear functional on $V$ and 
\[ 
\left\vert I_n(f) \right\vert = \left\vert \sum_{i=0}^{n} {w_{i}^{(n)}
    f\left(x_{i}^{(n)}\right)} \right\vert \leq \sum_{i=0}^{n} \left\vert
  w_{i}^{(n)} \right\vert \left\vert f\left(x_{i}^{(n)}\right)\right\vert
\leq {\left\Vert f \right\Vert}_{\infty} \sum_{i=0}^{n} \left\vert
  w_{i}^{(n)} \right\vert \, .\]
Thus $\left\Vert I_n \right\Vert = \sum_{i=0}^{n} \left\vert
  w_{i}^{(n)} \right\vert$, and so each $I_n$ is also a bounded linear
functional.

With these preliminaries, we can define the stability, consistency and
convergence of quadrature rules in exactly the same way as was done
in Definition~\ref{defn:ccs_linear}.

\begin{definition} 
  A quadrature rule $I_n$ is said to {\bf converge} to $I$ if $I_n(f)
  \rightarrow I(f)$ for every $f$ in $V$. It is {\bf consistent} if it
  converges on a dense subspace of $V$ and {\bf stable} if $\sup_{n}
  \left \Vert I_n \right\Vert < \infty$.
\end{definition}

\begin{remark}
  The motivation for the above definition of stability is the
  following.
\begin{align*}
  \left\vert I_n(f_1) - I_n(f_2)\right\vert &= \left \vert
    \sum_{i=0}^{n} {w_{i}^{(n)} f_{1}\left(x_{i}^{(n)}\right)} -
    \sum_{i=0}^{n} {w_{i}^{(n)}
      f_{2}\left(x_{i}^{(n)}\right)} \right\vert\\
  &= \left\vert \sum_{i=0}^{n} w_{i}^{(n)}\left(f_{1}\left(x_{i}^{(n)}
      \right) -
      f_{2}\left(x_{i}^ {(n)}\right)\right) \right\vert\\
  &\leq \sum_{i=0}^{n} \left\vert w_{i}^{(n)} \right\vert {\left\Vert
      f_1 -f_2 \right\Vert}_ {\infty}\, .
\end{align*}
Thus $\sup_{n} \sum_{i=0}^{n} \left\vert w_{i}^{(n)} \right\vert <
\infty $ gives discrete well-posedness, i.e., stability.
\end{remark}

Now we prove an equivalence theorem for quadrature whose statement and
proof is exactly the same same as the general linear equivalence
theorem proved earlier (Theorem~\ref{thm:equiv_linear}). We have
decided to use the natural number $n$ as a parameter here instead of
the real parameter $h$ of Theorem~\ref{thm:equiv_linear}, because this
is the natural setting for quadrature ($n+1$ is the number of nodes).

\begin{theorem}[\bf Equivalence Theorem for Quadrature]
  \label{thm:quadrature}
  A consistent quadrature rule is convergent if and only if it is
  stable.
\end{theorem}

\begin{proof}
  Suppose $I_n$ converges to $I$, i.e., $I_n(f) \rightarrow I(f)$ for
  every $f$ in $V$. This implies that for any given $f$ in $V$, the
  sequence $\left\{I_n(f)\right\}$ is bounded. Since each $I_n$ is a
  bounded linear functional, we can apply the uniform boundedness
  principle which gives us that $\sup_{n} \left\Vert I_n \right\Vert <
  \infty$ which is the definition of stability.

  Conversely assume stability and consistency. By definition of
  consistency then $I_n(f) \rightarrow I(f)$ for all $f$ in $V_0$
  where $V_0$ is a dense subspace of $V$. Stability means that
  $\sup_{n} \left\Vert I_n \right\Vert < \infty$.  By the density of
  $V_0$ in $V$, given $f \in V$ choose $f_{0} \in V_0 $ such
  that \[{\Vert f-f_0 \Vert}_{\infty} \leq \dfrac{\epsilon}{3 \max
    \{\Vert I \Vert, \sup_{n} \Vert I_n \Vert \}}\, .\]Hence
  \begin{align*}\Vert I (f) - I_n (f) \Vert &\leq \Vert I (f) - I
    (f_0) \Vert +
    \Vert I (f_0) - I_n (f_0) \Vert + \Vert I_n (f_0)  - I_n (f) \Vert \\
    & \leq \Vert I \Vert {\Vert f-f_0 \Vert}_{\infty} + \Vert I (f_0)
    - I_n (f_0) \Vert + \Vert I_n \Vert {\Vert f_0 -f \Vert}
    _{\infty}\\
    &\leq \Vert I \Vert {\Vert f-f_0 \Vert}_{\infty} + \Vert I (f_0) -
    I_n (f_0) \Vert + \sup_{n}
    \Vert I_n \Vert {\Vert f-f_0 \Vert}_{\infty}\\
    &\leq \dfrac{\epsilon}{3}+ \Vert I (f_0) - I_n (f_0) \Vert
    +\dfrac{\epsilon}{3} \, .
\end{align*}
By consistency, there exists $n_0 \in \mathbb{N}$ such that $\Vert I
(f_0) - I_n (f_0) \Vert \leq \epsilon/3$ for all $n \geq n_0$. Hence
$\Vert I (f) - I_n (f) \Vert \leq \epsilon$ for all $n \geq n_0$.
Therefore $I_n(f) \rightarrow I(f)$ for every $f$ in~$V$.
\end{proof}

Now we use the equivalence theorem to show convergence of Gaussian
quadrature rules and of the composite trapezoidal rule, followed by
the non convergence and instability of Newton-Cotes rules.

\begin{example}[Gaussian Quadrature]
  In Gaussian quadrature all the weights $w_{i}^{(n)} \geq 0$. Hence
  $\left\Vert I_n \right\Vert = \sum_{i=0}^{n} w_{i}^{(n)}$.  Gaussian
  quadrature rule $I_n$ is exact for all polynomials of degree less or
  equal to $2n-1$, i.e., $I_n(p) = \int_{a}^{b} p(x) dx$ for all
  polynomials $p$ of degree less than or equal to $2n-1$.  Since the
  space of such polynomials is dense in $V$, Gaussian quadrature is
  consistent. Moreover, $I_n(1) = \sum_{i=0}^{n} w_{i}^{(n)} =
  \int_{a}^{b}1 dx = b-a$. Thus $\sup_{n} \left\Vert I_n \right\Vert =
  \sup_{n} \sum_{i=0}^{n} w_{i}^{(n)} = b-a$. Therefore Gaussian
  quadrature is stable. Then by the above theorem, it is convergent.
\end{example}

\begin{example}[Composite Trapezoidal Rule]
  The composite trapezoidal rule has non negative weights. Moreover,
  it is exact for all piecewise linear polynomials which is a dense
  subspace of $V$. Hence by the argument as in the above example, we
  have stability and consistency and therefore convergence of the
  composite trapezoidal rule.
\end{example}

\begin{example}[Newton-Cotes]
  Define $I_n$ to be the Newton-Cotes quadrature rule. The nodes are
  are equally spaced in this quadrature rule and it integrates
  polynomials of a certain degree exactly. Thus $I_n$ is a consistent
  family by our definition. However it is not convergent.  An example
  continuous function for which $I_n (f)$ does not converge to $I(f)$
  is the Runge's function $f(x) = (1/1+25x^2)$ in the interval
  $[-1,1]$ (see \cite{SuMa2003}, page 208).  This function also
  appears in the polynomial interpolation section in
  Example~\ref{eg:runge_interp}. By Theorem~\ref{thm:quadrature}
  Newton-Cotes should also be unstable.  Indeed it is known that in
  Newton-Cotes rules some of the weights have negative sign and that
  this leads to instability, i.e., it is known that $\sup_{n}
  \left\Vert I_n \right \Vert = \sup_{n} \sum_{i=0}^{n} \vert
  w_{i}^{(n)} \vert = \infty$ (see page 350 of \cite{Heath2002}).
\end{example}

\subsection{Monte Carlo Integration}\label{subsec:mc}
For well behaved functions, i.e., functions with continuous
derivatives, the deterministic quadrature rule is very efficient at
least in one dimension. However, if the function fails to be well
behaved or in the case of multidimensional integrals, other techniques
can be competitive. In this section we will define convergence,
consistency and stability for Monte Carlo integration. For simplicity,
this is done for functions in $C[a,b]$. The notation and results from
probability theory that are needed were reviewed in
Section~\ref{sec:prelim}.

Let $V = (C[a,b], {\Vert \cdot \Vert}_{\infty})$ and for an $f \in V$,
let $I: V\rightarrow \mathbb{R}$ be defined as $I(f)= \int_{a}^{b}f$.
Let $ (W, {\Vert \cdot \Vert}_{{\infty}'})$ denote the space all of
bounded random variables defined on a probability space $(\Omega,
\Sigma, P)$, where ${\Vert X \Vert}_{{\infty}'} =
\operatorname{ess. sup_{\omega \in \Omega}} \vert X(\omega)\vert$.
Let $M_n : V \rightarrow W$ be a sequence of maps, not necessarily
bounded linear. For a specific Monte Carlo integration method see
Example~\ref{eg:sample_mean} below.  In that example the discrete
operators $M_n$ are bounded linear. We have chosen to state and prove
the equivalence theorem for Monte Carlo integration
(Theorem~\ref{thm:montecarlo}) without this assumption. This is done
to illustrate how the nonlinear case proof of an equivalence theorem
works in a probabilistic setting such as this.

We can look upon $I$ as a map from $V$ to $W$ by defining $I(f)$ to be
the constant random variable, i.e., $I(f): \Omega \rightarrow
\mathbb{R}$ is defined as $I(f)(\omega) = I(f)$ for all $\omega \in
\Omega$.

\begin{definition}
  A Monte Carlo integration is said to be {\bf convergent} if for any
  given $f \in V$,
  \[
  P \left[ \omega \in \Omega : \lim_{n \rightarrow \infty}
    M_n(f)(\omega) = I(f)(\omega) \right] =1 \, .
  \]
  It is {\bf consistent} if there is a dense subspace $V_0 $ of $V$
  such that for any $f \in V_0$, it is convergent. It is said to be
  {\bf stable at} $f_0 \in V$ if for any $\epsilon >0$ there exists a
  $\delta > 0$ such that for each $f \in V$ with ${\Vert f - f_0
    \Vert}_{\infty} \leq \delta$, there exists $n_0 \in \mathbb{N}$
  such that for all $n \geq n_0$,
  \[
  P \left[ \omega \in \Omega : \left \vert
      M_n(f)(\omega) - M_n(f_0)(\omega) \right \vert \leq \epsilon
  \right] =1\, .
  \] 
  It {\bf stable} if it is stable at every $f_0 \in V$.
\end{definition}

Now we state and prove an equivalence theorem for which the proof is
similar to Theorem~\ref{thm:equiv_nonlinear}, but with a probabilistic
flavor.

\begin{theorem}[\bf Equivalence Theorem for Monte Carlo
  Integration] \label{thm:montecarlo} A consistent Monte Carlo
  integration is convergent if and only if it is stable.
\end{theorem}

\begin{proof}
  Suppose it is convergent. To show stability, we need to find an $n_0
  \in \mathbb{N}$ and $\delta > 0$ such that for a given $f \in V$
  with ${\Vert f - f_0 \Vert}_{\infty} \leq \delta$, and for all $n
  \geq n_0$, we have $ P \left[ \omega : \left \vert M_n(f)(\omega) -
      M_n(f_0)(\omega) \right \vert \leq \epsilon \right] =1$.
  Outside a set of $P$-measure 0 in $\Omega$ and for all $n \in
  \mathbb{N}$, we have the following inequality
  \begin{align*}
    \left\vert M_n(f)(\omega)-M_n(f_0)(\omega) \right\vert
    &\leq \left\vert M_n(f)(\omega)- I(f)(\omega) \right\vert \\
    &+\left\vert I(f)(\omega) - I(f_0)(\omega) \right\vert +\left
      \vert I(f_0)(\omega)-M_n(f_0)(\omega) \right \vert \, .
  \end{align*}
  By convergence, there are $n_1$ and $n_2$ such that for all $n \geq
  n_1$, 
  \[
  P \left [ \omega : \left \vert M_n(f)(\omega) - I(f)(\omega) \right
    \vert \leq \epsilon/3 \right ] = 1 \, ,
  \]
  and for all $n \geq n_2$, 
  \[
  P \left [ \omega : \left \vert M_n(f_0)(\omega) - I(f_0)(\omega)
    \right \vert \leq \epsilon/3 \right ] = 1 \, .
  \]
  The integral operator is bounded, because if ${\Vert
    f-f_0\Vert}_{\infty} \leq \delta$, then $\vert I(f)-I(f_0) \vert
  \leq \delta(b-a)$.  Since $I(f)(\omega) = I(f)$ for all $\omega \in
  \Omega$ and by the boundedness of the integral operator, for an
  appropriately chosen $\delta > 0$, we have $\vert I(f)-I(f_0) \vert
  \leq \epsilon/3$. Hence
  \[
  P \left [ \omega : \left \vert I(f)(\omega) - I(f_0)(\omega)
    \right \vert = \left \vert I(f) - I(f_0) \right \vert \leq
    \epsilon/3 \right ] =1\, .
  \]
  Define the three sets 
  \begin{align*}
    \Omega_1 &= \{ \omega \in \Omega : 
    P \left [ \omega : \left \vert M_n(f)(\omega) - I(f)(\omega)
      \right \vert \leq \epsilon/3 \right ] =1 \}\\
    \Omega_2 &= \{ \omega \in \Omega :
    P \left [ \omega : \left \vert I(f)(\omega) - I(f_0)(\omega) \right
      \vert \leq \epsilon/3 \right ] =1\}\\
    \Omega_3 &= \{ \omega \in \Omega :
    P \left [ \omega : \left \vert I(f_0)(\omega) - M_n(f_0)(\omega)
      \right \vert \leq \epsilon/3 \right ] =1 \}\, .
  \end{align*}
  It is easy to see that $\Omega_2=\Omega$. Let $\Omega' =
  \bigcap_{i=1}^3 \Omega_i$. Let $n_0 = \max(n_1, n_2)$. For all $n
  \geq n_0$, $P(\Omega_i) =1$ for $1 \leq i \leq 3$.  Since $P$ is
  countably additive, measure of a countable union of sets of measure
  zero is zero. Therefore $P(\Omega')=1$. Thus for all $\omega \in
  \Omega'$, except for set of $P$-measure zero we have $\left \vert
    M_n(f)(\omega)-M_n(f_0)(\omega) \right\vert \leq \epsilon$ for all
  $n \geq n_0$. Hence $P \left[ \omega : \left \vert M_n(f)(\omega) -
      M_n(f_0)(\omega) \right \vert \leq \epsilon \right] =1$ for all
  $n \geq n_0$. Since $f_0$ was arbitrary the Monte Carlo integration
  is stable if it is convergent.

  Conversely, suppose we assume stability and consistency.  Therefore,
  given $f_0$ and $\epsilon > 0$, there exists $\delta > 0$, such that
  for each $f \in V$ with ${\Vert f - f_0 \Vert}_{\infty} \leq
  \delta$, there exists an $n_0 \in \mathbb{N}$ such that for all $n
  \geq n_0$,
  \[
  P \left[
    \omega \in \Omega : \left \vert M_n(f)(\omega) - M_n(f_0)(\omega)
    \right \vert \leq \epsilon \right] =1\, .
  \]
  By the density of $V_0$ in $V$, we can choose $g \in V_0$ such that
  ${\Vert g - f_0 \Vert}_{\infty} \leq \delta$.  By the boundedness of
  $I$, for an appropriately chosen $\delta$, we have ${\left \vert
      I(g)-I(f_0) \right \vert} \leq \epsilon/3$.  Again, outside a
  set of $P$-measure 0 in $\Omega$ and for all $n \in \mathbb{N}$, we
  have the following inequality
  \begin{align*}
    \left\vert M_n(f_0)(\omega)-I(f_0)(\omega) \right\vert
    &\leq \left\vert M_n(f_0)(\omega)- M_n(g)(\omega) \right\vert \\
    &+\left\vert M_n(g)(\omega) - I(g)(\omega) \right\vert +\left
      \vert I(g)(\omega)-I(f_0)(\omega) \right \vert \, .
  \end{align*}
  Stability at $f_0 $ means that there exists $\delta > 0 $ such that
  for each $g \in V_0$ with ${\Vert g - f_0 \Vert}_{\infty} \leq
  \delta$ there exists $n_1 \in \mathbb{N}$ with
  \[
  P \left [ \omega : \left \vert M_n(f_0)(\omega) - M_n(g)(\omega)
    \right \vert \leq \epsilon/3 \right ] =1 \}\, ,
  \]
  for all $n \geq n_1$.  By consistency, there is an $n_2$ such that
  for all $n \geq n_2$,
  \[
  P \left [ \omega : \left \vert M_n(g)(\omega) - I(g)(\omega) \right
    \vert \leq \epsilon/3 \right ] =1 \}\, .
  \]
  Moreover, 
  \[
  P \left [ \omega : \left \vert I(g)(\omega) - I(f_0)(\omega) \right
    \vert \leq \epsilon/3 \right ] =1\}\, .
  \]
  By an identical argument as in the previous half of the proof, on a
  subset $\Omega'$ of $\Omega$ with $P$-measure one, we have
  \[
  P \left [ \omega \in \Omega' : \left \vert M_n(f_0)(\omega) -
      I(f_0)(\omega) \right \vert \leq \epsilon \right ] =1 \}\, ,
  \]
  for all $n \geq n_0 = \max(n_1, n_2)$. Hence we have convergence at
  $f_0$. Since $f_0$ was arbitrary we have that stability implies
  convergence.
\end{proof}

As an example of numerical integration using Monte Carlo methods we
will examine the Sample-Mean Monte Carlo method and discuss its
convergence and stability. We prove both convergence and stability
separately. An alternative would have been to prove consistency
and one of the other two properties and infer the third from
Theorem~\ref{thm:montecarlo}.

\begin{example}[Sample-Mean Monte Carlo Method]
  \label{eg:sample_mean}
  Suppose we want to compute the approximate value of the integral $I=
  \int_{a}^{b} f(x) dx$. First, we choose any function $g \in C[a,b]$
  with the property that $g > 0$ and $\int_{a}^{b}g =1$.  Since $g \in
  C[a,b]$ and is positive, there exists $m \in \mathbb{R} $ such that $
  0 < m \leq g(x) $ in $[a,b]$.

  Then there exists some random variable $X$ with range in $[a,b]$,
  such that $g(x)$ is the pdf of $X$\cite{Tucker1967}.  Then, consider
  the random variable $Y= f(X)/g(X)$. Therefore,
  \[
  E(Y)= \int_{a}^{b}\dfrac{f(x)}{g(x)}g(x)dx = \int_{a}^{b}f =I\, .
  \]
  Now, choose a sequence of independent and identically distributed
  random variables $X=X_1$,$ X_2$, $\ldots$,$X_n$, $\ldots$ Since
  $f,g$ are continuous, they are Borel measurable. Since $g \neq 0$,
  $f/g$ is Borel measurable.  Hence the sequence of random variables
  $Y = f(X)/g(X)= f(X_1)/g(X_1)= Y_1, Y_2 = f(X_2)/g(X_2), \ldots,
  Y_n=f(X_n)/g(X_n), \ldots$ are independently and identically
  distributed by Theorem \ref{thm:independent}.  Hence for the chosen
  pdf $g$, we can define $M_{n}: V\rightarrow W$ as \[M_{n}(f) =
  \dfrac{1}{n}\sum _{i=1}^{n}\dfrac{f(X_i)}{g(X_i)}\, .\]By Theorem
  \ref{thm:lawof},
  \[
  P\left[\omega : \lim_{n \rightarrow \infty}\dfrac{1}{n}\sum
    _{i=1}^{n}Y_i(\omega) = I \right]=1\, .
  \]
  Thus we have the convergence of $M_{n}$. Since $M_{n}$ is linear in
  $f$, we can check for the boundedness of this linear operator
  \[
  \left \Vert M_{n} \right \Vert = \sup_{{\Vert f \Vert}_{\infty}=1}
  {\left \Vert \dfrac{1}{n}\sum _{i=1}^{n}\dfrac{f(X_i)}{g(X_i)}
    \right \Vert}_{{\infty}^{'}} \leq \dfrac{1}{n} \sum _{i=1}^{n}
  {\left \Vert \dfrac{f}{g} \right \Vert}_{\infty} \leq \dfrac{1}{m}
  \, .
  \]
  To get the estimate on the norm of $M_n$ we have used the fact that
  ${\Vert h(X)\Vert}_{\infty'} \leq {\Vert h \Vert}_{\infty}$, where
  $h = f/g$. This is true because
  \[ {\Vert h(X) \Vert}_{\infty'}= \operatorname{ess. sup}_{\omega \in
    \Omega} \vert h(X(\omega)) \vert = \inf \{M : P[\omega \in \Omega
  : \vert h(X(\omega))\vert \geq M ] =0 \} \, ,
  \]
  and hence the bound on $h(X)$ is controlled by the bound on $h$.

  The bound on $M_{n}$ is independent of $n$ and depends only on $g$.
  To exhibit stability we have to show that given $\epsilon > 0$ there
  exists $\delta > 0, n_0 \in \mathbb{N}$ such that for all $n \geq
  n_0$
  \[
  P \left[ \omega : \left \vert M_{n}(f)(\omega) - M_{n}(f_0)(\omega)
    \right \vert \leq \epsilon \right] =1\, ,
  \]
  whenever ${\Vert f - f_0 \Vert}_{\infty} \leq \delta$.  But outside
  a set of $P$-measure zero and for all $n \in \mathbb{N}$, we have
  \[
  \left \vert M_{n}(f)(\omega) - M_{n}(f_0)(\omega) \right \vert \leq
  {\left\Vert M_{n}(f) - M_{n}(f_0) \right\Vert}_{\infty'} \leq
  \left\Vert M_{n} \right\Vert {\Vert f- f_0 \Vert}_{\infty}\, .
  \]
  Therefore, for an appropriately chosen $\delta > 0 $ and for all $n
  \in \mathbb{N}$, we have
  \[
  P \left[ \omega : \left \vert M_{n}(f)(\omega) - M_{n}(f_0)(\omega)
    \right \vert \leq \epsilon \right] =1\, ,
  \]
  hence stability.

  The random variables $X_n$ are not necessarily unique for a given
  $g$. However, it does not matter as the Sample Mean Monte Carlo
  method is stable and convergent if we choose a positive pdf $g \in
  C[a,b]$.
\end{example}

\section{Numerical Differentiation}
\label{sec:differentiation}
If sufficiently differentiable functions are considered and a sum of
sup norms on the function and its derivatives is used as the norm then
smooth differentiation is bounded linear.  Numerical differentiation
can be posed as a parameterized collection of linear operators. If
they are assumed to be bounded as well then the definitions and
equivalence theorem from Section~\ref{sec:ccs} apply here. We show in
this section, and it is no surprise, that the equivalence theorem is
actually not needed for proving convergence or stability for the usual
finite difference formulas for the first derivative such as forward,
backward and central difference formulas. We also give the example of
the lowest order, 3 points finite difference formula for second
derivative, for which equivalence theorem is also not required. For
these and similar simple formulas, the proofs of stability and
convergence can be done directly and independently of each other and
are easy.

Thus there is no practical benefit of an equivalence theorem in the
context of such simple finite difference formulas for numerical
differentiation in one dimension. However the equivalence theorem
might be of benefit in proving stability or convergence for formulas
of high order accuracy for arbitrary derivatives on non equally spaced
grids such as the finite difference formulas in
\cite{Fornberg1988,Fornberg1996}.  The mesh size parameter $h$ appears
as $h^k$ in the denominator of these and similar formulas and
stability proofs might be tedious. Here $k$ is the order of the
derivative.  We do not discuss these formulas further in this paper.

For $k \geq 1$, define ${\Vert f \Vert}_{C^k} := \sum_{i=0}^{k} {\Vert
  f^{(i)} \Vert}_{\infty}$ where $f^{(i)}$ denotes the $i$-th
derivative of $f$. Let $V^k = (C^k [a,b], \Vert \cdot \Vert_{C^k})$
and let $W = (C[a,b], \Vert \cdot \Vert_{\infty})$. Consider
$D_{h}^{(k)}, D^{(k)} : V^k \rightarrow W$, where $D_{h}^{(k)}$ and
$D^{(k)}$ are the discrete and smooth differentiation operators
respectively. Here $h \in (0,1)$, is a parameter of the discrete
operators and is a measure of how far apart the points used in, say, a
finite difference formula are.  In the chosen norms, $D^{(k)}$ is 
a bounded linear operator.  In addition, assume also that $D_{h}^{(k)}$
are bounded linear operators.

The definitions of consistency, convergence and stability in 
Definition~\ref{defn:ccs_linear} can now be repeated in the
context of numerical differentiation.

\begin{definition}
  Numerical differentiation is said to be {\bf convergent} if
  \[
  \lim_{h\rightarrow 0}{\left\Vert \left(D_h^{(k)}-D^{(k)}\right)f\right
    \Vert}_{\infty} =0\, , 
  \]
  for every $f \in V^k$ and it is {\bf consistent} if it converges on
  a dense subspace of $V^k$. It is said to be {\bf stable} if $\Vert
  D_h^{(k)} \Vert \leq C$ where $C$ is independent of the
  parameter~$h$.
\end{definition}

\begin{theorem}[\bf Equivalence Theorem for Numerical
  Derivatives] \label{thm:differentiation} A consistent finite
  difference scheme is convergent if and only if stable.
\end{theorem}

\begin{proof}
  The proof is identical to the proof of
  Theorem~\ref{thm:equiv_linear}, with $T_h$ replaced by $D_h^{(k)}$
  and $T$ by $D^{(k)}$.
\end{proof}

\begin{example}[Forward Difference]
  As an example we now consider the basic forward difference
  approximation to the first derivative which we will show to be
  stable and convergent. The equivalence theorem is not required in
  this case although using it would reduce the work required in
  proving stability and convergence. It is easy enough to prove
  convergence and stability separately in an elementary way. Consider
  a point $x \in [a,b-h]$.  Let $h \in (0,1)$ such that $x+h < b$.  We
  will show that the forward difference approximation to the first
  derivative
  \[
  D_h^{(1)}f(x) := \dfrac{f(x+h)-f(x)}{h}\, , 
  \]
  is both convergent and stable. First note that
  \begin{align*}
    \left\Vert D_h^{(1)}\right\Vert=\sup_{{\left\Vert f
        \right\Vert}_{C^1}=1} {\left\Vert D_h^{(1)}f
      \right\Vert}_{\infty} &=\sup_{{\left\Vert f \right\Vert}_{C^1}=1}
    \sup_{x \in [a,b-h]}\left\vert \dfrac{f(x+h)-f(x)}{h}
    \right\vert \\
    & = \sup_{{\left\Vert f \right\Vert}_{C^1}=1} \sup_{x \in [a,b-h]}
    \left\vert f'(x + \theta h) \right\vert \leq 1 \, ,
  \end{align*}
  for some $0 < \theta <1$. Thus $D_h^{(1)}$ is stable. Moreover, by
  Mean Value Theorem we have
  \begin{align*}
    \lim_{h\rightarrow 0}{\left\Vert D_h^{(1)}f-D^{(1)}f\right\Vert}_{\infty} 
    &= \lim_{h\rightarrow 0} \sup_{x \in [a,b-h]}
    \left\vert \dfrac{f(x+h)-f(x)}{h}-f'(x) \right\vert\\
    &=\lim_{h\rightarrow 0} \sup_{x \in [a,b-h]}
    \left\vert f'(x + \theta h)-f'(x) \right\vert =0 \, ,
  \end{align*}
  which means that $D_h^{(1)}$ converges to $D^{(1)}$.  The
  convergence and stability of the other commonly used finite
  difference schemes like backward difference and central difference
  schemes can be similarly proved.
\end{example}

The stability proved above means that a slight perturbation of $f \in
V^1$ does not drastically change the computed numerical
derivative. Note however that if the $C^{1}$ norm is replaced by the
sup norm, i.e., if the space $V^1$ is replaced by $(C^1[a,b], {\Vert \cdot
  \Vert}_{\infty})$, then the finite difference schemes above are
highly sensitive to small perturbations, leading to instability. This
can be seen in the following example.

\begin{example}[Unboundedness in sup Norm]\label{eg:unbounded_deriv}
  Define a sequence $g_h(x)= \sin (2\pi x/h)$ in $\left(C^{1}[0,1],
    \Vert \cdot \Vert_{\infty}\right)$.  Note that $\left\Vert g_h
  \right\Vert_{\infty}=1$.  However,
  \[
  {\left\Vert D_h^{(1)}(g_h) \right\Vert}_{\infty} = 
  \sup_{x \in [0,1-h]}\left\vert\dfrac{\sin (2\pi(x+h)/h)- \sin (2\pi x/h)}
    {h}\right\vert\, .
  \] 
  By Mean Value Theorem,
  \[
  \dfrac{\sin (2\pi (x+h)/h)- \sin (2\pi x/h)}{h} = \dfrac{2\pi}{h}\cos
  (2\pi c/h) \, ,
  \]
  for some $c \in (0,1-h)$ and so
  \[
  {\left\Vert D_h^{(1)}(g_h) \right\Vert}_{\infty} 
  = \sup_{c \in (0,1-h)}\left\vert \dfrac{2\pi}{h}\cos (2\pi c / h)
  \right\vert  =\dfrac{2\pi}{h}  \, .  
  \]
  This implies that ${\Vert D_h^{(1)} \Vert}_{\infty} > 1/h$ for all $
  h \in (0,1)$.
\end{example}

\begin{example}[Second Derivative]
  We now discuss the convergence and stability of the most basic
  finite difference operator for second derivative. We will show that
  the finite difference approximation for second derivative
  \[
  D_h^{(2)}f(x) := \dfrac{f(x+h)-2f(x)+f(x-h)}{h^2}\, , 
  \]
  is both convergent and stable. The operators $D_h^{(2)}$ are
  consistent and so it is actually enough to prove just stability or
  convergence due to Theorem~\ref{thm:differentiation}.  But as in the
  first derivative case, we will prove stability and convergence
  separately, since both are easy to prove. To prove stability, note
  that
  \[
  \left\Vert D_h^{(2)} \right\Vert =\sup_{{\left\Vert f \right\Vert}_{C^2}=1}
  {\left\Vert D_h^{(2)}f \right\Vert}_{\infty}
  =\sup_{{\left\Vert f \right\Vert}_{C^2}=1}\sup_{x \in [a+h,b-h]}
  \left\vert \dfrac{f(x+h)-2f(x)+f(x-h)}{h^2} \right\vert \, .
  \]
  By two applications of the Mean Value Theorem, we have that for some
  $0 < \alpha, \beta <1$, the above is equal to
  \[
  \sup_{{\left\Vert f \right\Vert}_{C^2}=1}\sup_{x \in [a+h,b-h]}
  \left\vert \dfrac{f'(x+\alpha h)-f'(x- \beta h)}{h} \right\vert
  \leq \sup_{{\left\Vert f \right\Vert}_{C^2}=1} 
  {\left\Vert f'' \right\Vert}_{\infty} \leq 1\, . 
  \]
  Hence the $D_h^{(2)}$ is stable. Moreover, once again by the Mean Value
  Theorem, for some $0<\alpha, \beta, \gamma <1$ we have
  \begin{align*}
    \lim_{h\rightarrow 0}{\left\Vert D_h^{(2)}f-D^{(2)}f \right\Vert}_{\infty}&=
    \lim_{h\rightarrow 0} \sup_{x\in [a+h,b-h]}
    \left\vert\dfrac{f(x+h)-2f(x)+f(x-h)}{h^2}-f''(x) \right\vert\\
    &=\lim_{h\rightarrow 0} \sup_{x\in [a+h,b-h]}\left\vert
      \dfrac{f'(x+ \alpha h)-f'(x- \beta h)}{h}-f''(x) \right\vert\\
    &=\lim_{h\rightarrow 0} \sup_{x\in [a+h,b-h]}
    \left\vert f''(x-\beta h + \gamma h (\alpha + \beta))-
      f''(x)\right\vert=0 \, , 
  \end{align*}
  which means convergence.  Again, if the $C^2$ norm is replaced by
  the $C^1$ norm or the sup norm, the operator $D_h^{(2)}$ fails to be
  stable. Examples similar to Example~\ref{eg:unbounded_deriv} above
  can be constructed to exhibit the unboundedness of the operator
  $D_h^{(2)}$.
\end{example}

\section{Polynomial Interpolation}
\label{sec:interp}

The interpolation problem is to find a function that takes on
prescribed values at specified points. In one dimension the data is
given as $(x_i, y_i)$ for $i=0, \ldots,n$, with $x_0 < x_1 < \cdots <
x_n$, and we look for a function $f : \mathbb{R} \rightarrow
\mathbb{R}$ called interpolating function such that $f(x_i)= y_i$ for
all $i$.  In order to define the notions of consistency, convergence
and stability it is convenient to pose the interpolation problem as
interpolation of continuous functions. For any given data which is to
be interpolated, there is the obvious unique piecewise linear
continuous function that interpolates the data. Stability with respect
to changes in the given values or the locations of the values are both
captured by stability with respect to changes in the piecewise linear
continuous function. We will only address interpolation of continuous
functions by polynomials.  It is easy to show the existence and
uniqueness of the interpolant polynomial in the case of one
dimensional interpolation, i.e., when the dimension of the domain of
the function is one \cite{Heath2002, Walsh1965}.

Let $\{p_n\}$ be a sequence of polynomials of degree at most $n$,
interpolating a function $f$ in $V = \left(C[a,b], {\left\Vert \cdot
    \right \Vert}_{\infty}\right)$ on a set of interpolation points
$x_i^{(n)}$, where $a \leq x_0^{(n)} < x_1^{(n)} < \cdots < x_n^{(n)}
\leq b$.  We can define polynomial interpolation as a map $P_n : V
\rightarrow \mathbb{P}_n \subset V$, where $\mathbb{P}_n$ is the space
of polynomials of degree at most $n$ and $P_n f $ is defined as the
unique polynomial $p \in \mathbb{P}_n$ which interpolates $f$ at the
given nodes. A basic fact about error in polynomial interpolation is
given by the following lemma.

\begin{lemma}\label{lem:interp_error} 
  Let $f \in C^{(n+1)}[a,b]$, and let $a \leq x_0 < x_1 < \cdots < x_n
  \leq b$ be a partition of $[a,b]$. Then for all $x \in (a,b)$, there
  exists $\xi \in (\min\{x_0,x\}, \max\{x_n,x\})$ such that
  \[
  f(x) - (P_n f)(x) = \frac{f^{(n+1)}(\xi)}{(n+1) !}
  \prod_{i=0}^{n}(x-x_i)\, .
  \]
\end{lemma}

\begin{proof}
  The proof requires repeated use of Rolle's Theorem. See page 119 of
  \cite{AtHa2005} for details.
\end{proof}

The following Lemma~\ref{lem:interp_is_linear} shows that the $P_n$
operators are linear, and Lemma~\ref{lem:interp_is_bounded} shows that
they are bounded. The corresponding smooth operator is the identity
map which is bounded linear.  Thus the linear case definitions of
consistency, convergence and stability given in
Definition~\ref{defn:ccs_linear} apply here as well.

\begin{lemma}[Linearity]\label{lem:interp_is_linear}
  The interpolation operator $P_n : V \rightarrow \mathbb{P}_n \subset
  V$ is linear.
\end{lemma}
\begin{proof} 
  Let $P_n(f) = p_1$ and $P_n(g)=p_2$, where $f,g \in V$ and $p_1,
  p_2$ are the unique polynomials in $\mathbb{P}_n$ such that $f(x_i)
  = p_1(x_i)$ and $g(x_i) = p_2(x_i)$ for $0 \leq i \leq n$. Hence
  $(f+g)(x_i) = (p_1+p_2)(x_i)$. By uniqueness of polynomial
  interpolation $P_n(f+g) = p_1+p_2$. Hence $P_n(f +g) =
  P_nf+P_ng$. Similarly using uniqueness we can show that $P_n(cf) =
  cP_n(f)$, for any $c \in \mathbb{R}$.
\end{proof}

\begin{lemma}[Boundedness]\label{lem:interp_is_bounded}
  Each interpolation operator $P_n$ is bounded, i.e., $ \Vert P_n
  \Vert = \sup_{{\Vert f \Vert}_{\infty}=1} {\Vert P_n(f)
    \Vert}_{\infty}\leq K(n)$ where $K(n) > 0$.
\end{lemma}
\begin{proof} 
  By Lemma~\ref{lem:interp_error}, $P_n (f) = f$, where $f$ is any
  polynomial of degree $\leq n$. Hence $P_n$ maps onto $\mathbb{P}_n$.
  Let $f_m$ be a sequence in $ P_n^{-1}(0)$ converging to $f$ in $
  V$. Since $f_m(x_i)=0$ for all $0 \leq i \leq n$ and $\lim_{m
    \rightarrow \infty}f_m(x) =f(x)$ for all $x \in [a,b]$, we have
  $f(x_i)= 0$ for all $0 \leq i \leq n$.  Therefore $f \in
  P_n^{-1}(0)$.  Hence the kernel of $P_n$ is closed.  Therefore, by
  Lemma~\ref{lem:closed_kernel}, $P_n$ is continuous and hence
  bounded.
\end{proof}

\begin{definition} 
  Interpolation is {\bf convergent} if for any given $f \in V$ the
  sequence of interpolant polynomials $P_n f$ converges to $f$ in
  $V$. It is {\bf consistent} if there is a dense subspace $V_0$ of
  $V$ such that for any $f$ in $V_0$, the sequence of interpolant
  polynomials $P_n f$ converges to $f$ in $V$.  Interpolation is said
  to be {\bf stable} if $\sup_n\Vert P_n\Vert < \infty$.
\end{definition}

If the function being interpolated is a polynomial then the
interpolant is exact. This follows from Lemma~\ref{lem:interp_error}.
Since polynomials are dense in $V$ (Weierstrass Approximation Theorem)
this implies that interpolation of continuous functions by polynomials
is consistent.

We now prove an equivalence theorem for polynomial interpolation. In
the notation of Theorem~\ref{thm:equiv_linear} the discrete operators
$T_n$ of the theorem correspond to the polynomial interpolation
operator $P_n$ and the smooth operator $T$ corresponds to the identity
operator.

\begin{theorem}[\bf Equivalence Theorem for Interpolation]\label{thm:interp}
  A consistent polynomial interpolation is convergent if and only if
  it is stable.
\end{theorem}

\begin{proof}
  The proof is identical to the one given for the general linear
  equivalence theorem (Theorem~\ref{thm:equiv_linear}) with
  discretization parameter $h$ replaced by $n$. Thus the proof of the
  equivalence theorem for quadrature (Theorem~\ref{thm:quadrature})
  can be used here with appropriate modifications.
\end{proof}

\begin{example}[Runge's Function] \label{eg:runge_interp}
  It is well known that for the Runge's function $f(x)= 1/(1+25x^2)$
  in the interval $[-1,1]$ the polynomial interpolants do not converge
  if the interpolation points are uniformly spaced.  See for example
  \cite{IsKe1966}.  As we noted earlier, interpolation of continuous
  functions by polynomials is consistent. However, the method is not
  stable if equispaced points are used for interpolation. In light of
  the equivalence theorem above (Theorem~\ref{thm:interp})
  this agrees with the lack of convergence. To see the lack of
  stability, let $p(x)$ be any polynomial which is arbitrarily close
  to the Runge's function. Existence of $p(x)$ is guaranteed by the
  Weierstrass approximation Theorem. Now, the interpolating
  polynomials for $p(x)$ are exact as we increase the number of
  interpolating points. However, the interpolating polynomial for the
  Runge's function is not close to $p(x)$. This shows that polynomial
  interpolation of continuous functions using equispaced points is not
  discretely well-posedness, i.e., it is not stable. 
\end{example}

Interpolating real valued functions whose domain has dimension greater
than one is a much harder problem, unlike the one dimensional case
where polynomial interpolant always exists and is unique. For example,
in 2 variable interpolation, if data is specified at 3 points that are
collinear, then there are infinitely many linear interpolants.  The
geometry of the point locations becomes an important factor in
existence and uniqueness of the interpolant. For a result of this type
in the 2 variable case see, for example, \cite{Phillips2003}. We need
the existence and uniqueness of interpolant so that there is a well
defined operator $P_n$.

There does exist a generalization, called Ciarlet's error formula
\cite{CiWa1971,CiRa1972,GaSa2000}, of Lemma \ref{lem:interp_error} to
the multivariable case. Moreover, multivariate polynomials are dense
in the sup norm, in the space of multivariate continuous functions on
compact subsets of $\mathbb{R}^n$ (Stone-Weierstrass Theorem)
\cite{Davis1975}. Thus when multivariate polynomial interpolation does
exist and is unique, we get consistency as before. Consistency,
stability and convergence can be defined as described in the
univariate case and thus an equivalence theorem exists for the
multivariate polynomial interpolation of continuous functions. Of
course the conditions for stability and convergence are more
complicated in the multivariate case and we do not address those here.

In contrast with interpolation, approximation methods do not require
agreement with data at specific points. For example one might seek a
polynomial $p$ such that $\min_{p \in \mathbb{P}_n} {\Vert f-
  p\Vert}_{\infty}$ is attained. In this particular case one can show
existence and uniqueness. See for instance Theorem 7.5.6 in
\cite{Davis1975}.  One can define the notion of convergence,
consistency and stability for approximation problems in an identical
fashion as for interpolation and prove a theorem like
Theorem~\ref{thm:interp}.  In the particular case considered
above, the approximation always converges (Theorem 7.6.1 of
\cite{Davis1975}) and hence is stable.

\section{Conclusions and Future Work}
We have shown that equivalence theorems can be proved in a general
setting in numerical analysis. As in the classical Lax-Richtmyer
equivalence theorem in PDEs, these theorems state that a consistent
method is convergent if and only if it is stable.  The notion of
stability we used was that of discrete well posedness and we defined
consistency to be convergence on a dense subspace.  The
discretizations of the smooth problems were considered to be linear or
nonlinear operators on normed linear spaces and depending on some
parameter which measures the discretization size. We showed that our
general equivalence theorems require basic tools like uniform
boundedness principle, triangle inequality, and density
arguments. Convergence implying stability was always obtained
independent of consistency.

In this paper, we have studied stability with respect to perturbation
of input data of the discrete operators. However, one could also
study stability with respect to locations of points, shape of domain
etc. One can also investigate equivalence theorems in other numerical
analysis contexts. Some examples are, multidimensional quadrature
rules, multidimensional Monte Carlo integration, optimization,
eigenvalue problems in PDEs, formulas for numerical differentiation of
any order and accuracy, such as those given in \cite{Fornberg1996}.

We defined consistency as convergence on a dense subspace. In our examples,
many times the dense subspace was polynomials in continuous functions.
It may be worthwhile to study if and how the choice of particular dense
subspace affects the theory and applications of equivalence theorems.

\bibliographystyle{amsplain}
\bibliography{lax}

\end{document}